\documentclass[a4paper,12pt]{article} 

\usepackage{color}

\usepackage{amsfonts, amsmath, amsthm, amssymb}
\usepackage[T1]{fontenc}
\usepackage[cp1250]{inputenc}
\usepackage{xcolor}
\usepackage{graphicx}
\usepackage{amssymb}
\usepackage{amsmath}
\usepackage{mathptmx}
\usepackage{helvet}
\usepackage{courier}
\usepackage{txfonts}
\usepackage{tikz} 
\usetikzlibrary{arrows}
\usepackage{type1cm}

\usepackage{verbatim}

\usepackage{graphicx}
\usepackage{epsfig,amscd,amssymb,amsxtra,amsmath,amsthm}
\usepackage{type1cm}
\usepackage[T1]{fontenc}
\usepackage{graphics}
\usepackage[mathscr]{eucal}
\usepackage[all]{xy}
\usepackage{amsmath,amscd}

\newtheorem{theorem}{Theorem}[section]

\newtheorem{definition}[theorem]{Definition}
\newtheorem{lemma}[theorem]{Lemma}

\newtheorem{example}[theorem]{Example}

\newtheorem{corollary}[theorem]{Corollary}
\newtheorem{problem}[theorem]{Problem}
\newtheorem{observation}[theorem]{Observation}

\newcommand{\Cl}  {\mathop{\rm Cl}\nolimits}
\newcommand{\Int}  {\mathop{\rm Int}\nolimits}

\begin{document}

\def\joinrel{\mkern-3mu}
\newcommand{\varproj}{\displaystyle \lim_{\multimapinv\joinrel-\joinrel-}}

\title{Retract or Not: A Tale of Two Fans}
\author{Iztok Bani\v c, Goran Erceg, Sina Greenwood, and Judy Kennedy}
\date{}

\maketitle

\begin{abstract}
Let $X$ be a Lelek fan or a Cantor fan and let $Y$ be a Lelek fan or a Cantor fan. In this paper, we study embeddings \( f: X \to Y \) that admit retractions from \( Y \) onto \( f(X) \). In 1989, W.~J.~Charatonik and J.~J.~Charatonik proved that if \( X \) is a Lelek fan and \( Y \) is a Cantor fan, then no embedding \( f \) of \( X \) into \( Y \) admits a retraction from \( Y \) onto \( f(X) \). They also showed that if both \( X \) and \( Y \) are Cantor fans, then every embedding \( f \) of \( X \) into \( Y \) admits such a retraction.

In this paper, we address the two remaining cases. First, we consider the situation where \( X \) is a Cantor fan and \( Y \) is a Lelek fan. We prove that in this case, every embedding \( f \) of \( X \) into \( Y \) admits a retraction from \( Y \) onto \( f(X) \). Second, we examine the case where both \( X \) and \( Y \) are Lelek fans. Here, we show that there exist embeddings \( f \) that do admit a retraction from \( Y \) onto \( f(X) \), as well as embeddings that do not. For this latter case, we also identify additional properties of embeddings that ensure the existence of a retraction from \( Y \) onto \( f(X) \).
\end{abstract}
\-

\section{Introduction}
 The study of retractions in topological spaces, particularly in continuum theory, has a long and rich history. In this paper, we focus on the question of whether certain embeddings between well-known dendroidal continua admit retractions. More specifically, we consider the Lelek fan and the Cantor fan - two classical examples of smooth, uniquely arcwise connected continua with intricate topological structures.

Let \( X \) and \( Y \) be continua, and let \( f: X \to Y \) be an embedding. A natural question arises: does \( f(X) \) admit a retraction from \( Y \)? This question becomes particularly interesting when \( X \) and \( Y \) are either Lelek fans or Cantor fans.

In 1989, W.~J.~Charatonik and J.~J.~Charatonik investigated this problem and established two foundational results:
\begin{itemize}
    \item If \( X \) is a Lelek fan and \( Y \) is a Cantor fan, then no embedding \( f: X \to Y \) admits a retraction from \( Y \) onto \( f(X) \).
    \item If both \( X \) and \( Y \) are Cantor fans, then every embedding \( f: X \to Y \) admits a retraction from \( Y \) onto \( f(X) \).
\end{itemize}

These results leave open the two remaining cases: when \( X \) is a Cantor fan and \( Y \) is a Lelek fan, and when both \( X \) and \( Y \) are Lelek fans.
The aim of this paper is to address these two cases. We prove that in the first case, every embedding \( f: X \to Y \) admits a retraction from \( Y \) onto \( f(X) \). This is somewhat surprising, given the failure of retractions in the reverse situation, where the Lelek fan is embedded into the Cantor fan.
In the second case, where both \( X \) and \( Y \) are Lelek fans, the situation is more subtle. We show that there exist embeddings \( f: X \to Y \) that admit retractions, as well as embeddings that do not. We further investigate this case by identifying specific conditions on the embedding \( f \) under which a retraction from \( Y \) onto \( f(X) \) does exist.

This paper is organized as follows. In Section \ref{s1}, we recall the definitions and relevant properties of the Lelek fan and the Cantor fan. Section \ref{s2} contains the proof of the positive result for embeddings of Cantor fans into Lelek fans. In Section \ref{s3}, we analyze embeddings of Lelek fans into Lelek fans and present examples, counterexamples, and sufficient conditions for the existence of retractions. In Section \ref{s4}, we introduce and study simple and semi-simple retractions from Lelek fans into Lelek fans, characterize when such retractions factor through the corresponding fences, and conclude with an open problem.

\section{Definitions}\label{s1}
In this section, we recall the basic definitions and fundamental properties of the Lelek fan and the Cantor fan that are relevant to our study. We also summarize known results concerning retractions induced by embeddings between these continua. Our aim is to provide a self-contained background that will serve as the foundation for the results presented in the following sections. Unless stated otherwise, all spaces considered are compact and metrizable, and all maps are assumed to be continuous.
\begin{definition}
 \emph{A continuum} is a non-empty compact connected metric space.  \emph{A subcontinuum} is a subspace of a continuum, which is itself a continuum.
 \end{definition}
 
 \begin{definition}
Let $X$ be a continuum. 
\begin{enumerate}
\item The continuum $X$ is \emph{unicoherent} if for any subcontinua $A$ and $B$ of $X$ such that $X=A\cup B$,  the compactum $A\cap B$ is connected. 
\item The continuum $X$ is \emph{hereditarily unicoherent } provided that each of its subcontinua is unicoherent.
\item The continuum $X$ is a \emph{dendroid} if it is an arcwise connected, hereditarily unicoherent continuum.
\item If $X$ is homeomorphic to $[0,1]$, then $X$ is \emph{an arc}.   
\item Let $X$ be an arc. A point $x\in X$ is called \emph{an end point of $X$} if  there is a homeomorphism $\varphi:[0,1]\rightarrow X$ such that $\varphi(0)=x$.
\item Let $X$ be a dendroid.  A point $x\in X$ is called an \emph{end point of $X$} if for  every arc $A$ in $X$ that contains $x$, $x$ is an end point of $A$.  The set of all end points of $X$ is denoted by $E(X)$. 
\item The continuum $X$ is \emph{a simple triod} if it is homeomorphic to $([-1,1]\times \{0\})\cup (\{0\}\times [0,1])$.
\item Let $X$ be a simple triod. A point $x\in X$ is called \emph{the top point} or, briefly, \emph{the top of $X$} if  there is a homeomorphism $\varphi:([-1,1]\times \{0\})\cup (\{0\}\times [0,1])\rightarrow X$ such that $\varphi(0,0)=x$.
\item Let $X$ be a dendroid.  A point $x\in X$ is called \emph{a ramification point of $X$}, if there is a simple triod $T$ in $X$ with top $x$.  The set of all ramification points of $X$ is denoted by $R(X)$. 
\item The continuum $X$ is \emph{a  fan} if it is a dendroid with at most one ramification point $v$, which is called \emph{the top of the fan $X$} (if it exists).
\item Let $X$ be a fan.   For all points $x$ and $y$ in $X$, we define  \emph{$[x,y]$} to be the arc in $X$ with end points $x$ and $y$, if $x\neq y$. If $x=y$, then we define $[x,y]=\{x\}$.
\item Let $X$ be a fan with top $v$. For each end point $e\in E(X)$ of the fan $X$, we call the arc in $X$ from $v$ to $e$, \emph{a leg of $X$}. The set of all legs of $X$ is denoted by $\mathcal L(X)$.
\item Let $X$ be a fan with top $v$. We say that the fan $X$ is \emph{smooth} if for any $x\in X$ and for any sequence $(x_n)$ of points in $X$,
$$
\lim_{n\to \infty}x_n=x \Longrightarrow \lim_{n\to \infty}[v,x_n]=[v,x].
$$ 
\item Let $X$ be a fan.  We say that $X$ is \emph{a Cantor fan} if $X$ is homeomorphic to the continuum $\bigcup_{c\in C}A_c$, where $C\subseteq [0,1]$ is the Cantor middle third set and for each $c\in C$, $A_c$ is the  {convex} segment in the plane from $(\frac{1}{2},0)$ to $(c,1)$.
\item Let $X$ be a fan.  We say that $X$ is \emph{a Lelek fan} if it is smooth and $\Cl(E(X))=X$. See Figure \ref{figure2}.
\begin{figure}[h!]
	\centering
		\includegraphics[width=20em]{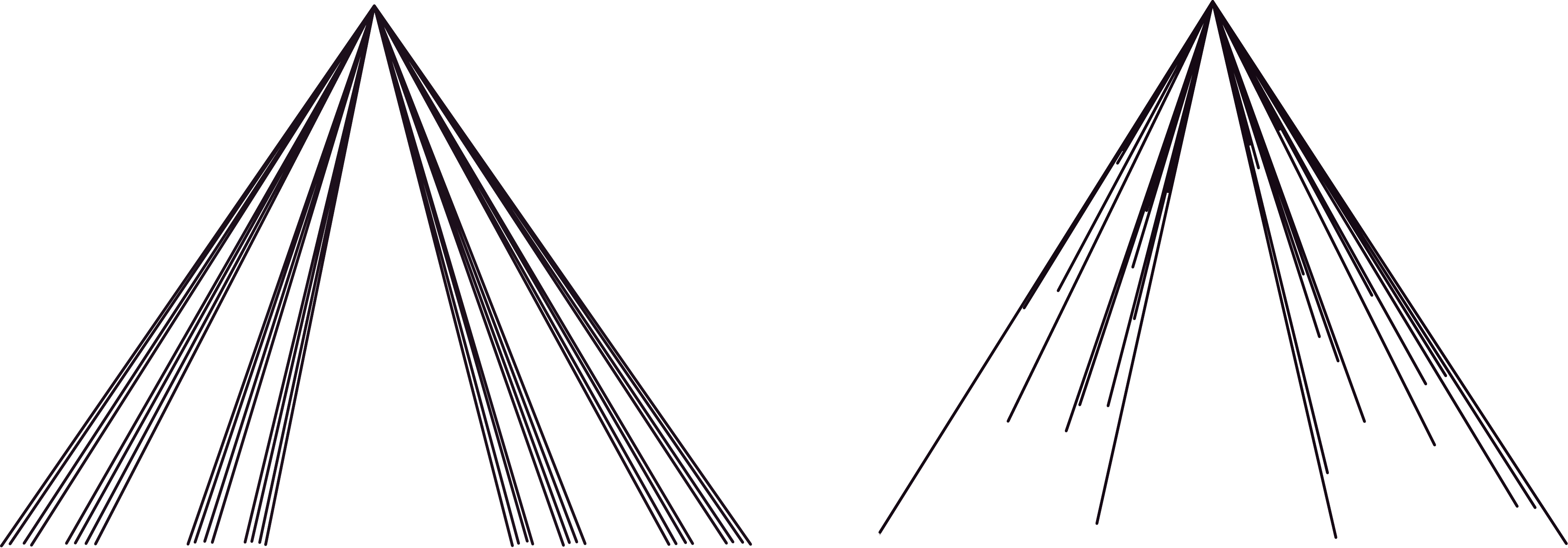}
	\caption{A Cantor fan and a Lelek fan}
	\label{figure2}
\end{figure}  
\end{enumerate}
\end{definition}
Note that a Lelek fan was constructed by A.~Lelek in \cite{lelek}.  An interesting  property of the Lelek fan $L$ is the fact that the set of its end points is a dense one-dimensional set in $L$. It is also unique, i.e., any two Lelek fans are homeomorphic, for the proofs see \cite{charatonik} and \cite{oversteegen}.
\begin{definition}
	Let $X$ and $Y$ be metric spaces and let $f:X\rightarrow Y$ be a continuous function. We say that $f$ is \emph{an embedding of $X$ into $Y$} if there is a homeomorphism $\varphi:X\rightarrow f(X)$ such that for each $x\in X$, $\varphi(x)=f(x)$.
\end{definition}
\begin{definition}
	Let $X$ and $Y$ be metric spaces. We say that $X$ is \emph{embeddable into $Y$} if there is an embedding of $X$ into $Y$.
\end{definition}
It easily follows from Lelek's construction \cite{lelek} that the Lelek fan is embeddable into the Cantor fan.  However, it is  not that obvious that the Cantor fan is embeddable into the Lelek fan.  {One can easily construct an embedding of the Cantor fan into the Lelek fan by using
\begin{enumerate}
\item the well-known result  from \cite{dijkstra}  by J. ~J.~ Dijkstra and J. ~Mill that a space is almost zero-dimensional \footnote{A space is called almost zero-dimensional if every point of the space has a neighbourhood basis consisting of C-sets of the space, where a subset $A$ of a space $X$ is called a C-set in $X$ if $A$ can be written as an intersection of clopen subsets of $X$; see [7] for more details.} if and only if it is embeddable into the complete Erd\"os space, and 
\item the well-known result from \cite{kawamura}  by K. ~Kawamura, L. ~G.~Oversteegen, and E.~D.~Tymchatyn that the set of end points of the Lelek fan is homeomorphic to the complete Erd\"os space.  
\end{enumerate} 
First, embed the Cantor set into the set of end points of the Lelek fan and then,   the subcontinuum  of the Lelek fan that is induced by the embedded Cantor set,  is a Cantor fan (among other things, this was already noted by G.~Basso and R.~Camerlo in \cite{basso},  where another similar result is obtained).} In \cite{EK}, another short and elegant proof that the Lelek fan can be embedded into the Cantor fan, is presented. In Theorem \ref{union}, we show not only that Lelek fans contain Cantor fans, we show that they are unions of Cantor fans.
\begin{theorem}\label{union}
Let $L$	be a Lelek fan with top $v$. Then for each $A\in \mathcal L(L)$ there is a Cantor fan $C_A$ in $L$ such that for all $A_1,A_2\in \mathcal L(L)$,
$$
A_1\neq A_2 ~~~  \Longrightarrow  ~~~ C_{A_1}\cap C_{A_2}=\{v\}  
$$
and
$$
\bigcup_{A\in \mathcal L(L)}C_A=L.
$$ 
\end{theorem} 
\begin{proof}
	Let $L'$ be a Lelek fan in $\mathbb R^2$, let $f:L\rightarrow L'$ be a homeomorphism, and let $C$ be the Cantor middle third set. We define the relation $\sim$ on $C\times L'$ by
	$$
	(c_1,x_1)\sim (c_2,x_2) ~~~ \Longrightarrow ~~~  (c_1,x_1)=(c_2,x_2) \textup{ or }  x_1=x_2=f(v).
	$$ 
	Note that $\sim$ is an equivalence relation on $C\times L'$ and that the quotient $(C\times L')/_{\sim}$ is a Lelek fan with top $\mathbf v=[(c,v)]$ for some $c\in C$. We leave the details to the reader. Let $F=(C\times L')/_{\sim}$ and let $\varphi:F\rightarrow L$ be a homeomorphism. Note that for each $A\in \mathcal L(L)$, $(C\times f(A))/_{\sim}$ is a Cantor fan in $F$ and that for all $A_1,A_2\in \mathcal L(L)$,
$$
A_1\neq A_2 ~~~  \Longrightarrow  ~~~ \left((C\times f(A_1))/_{\sim}\right)\cap \left((C\times f(A_2))/_{\sim}\right)=\{\mathbf v\}  
$$
and
$$
\bigcup_{A\in \mathcal L(L)}\left((C\times f(A))/_{\sim}\right)=F.
$$ 
   For each $A\in \mathcal L(L)$, let $C_A=\varphi((C\times f(A))/_{\sim})$. Then for each $A\in \mathcal L(L)$,  $C_A$ is a Cantor fan in $L$ such that for all $A_1,A_2\in \mathcal L(L)$,
$$
A_1\neq A_2 ~~~  \Longrightarrow  ~~~ C_{A_1}\cap C_{A_2}=\{v\}  
$$
and
$$
\bigcup_{A\in \mathcal L(L)}C_A=L.
$$
\end{proof}
\begin{definition}
	Let $X$ be a metric space, let $Y$ be a subspace of $X$, and let $r:X\rightarrow Y$ be a function. If $r$ is continuous such that for each $y\in Y$, $r(y)=y$, then we say that $r$ is \emph{a retraction from $X$ onto $Y$}. 
\end{definition}
\begin{definition}
	Let $X$ and $Y$ be metric spaces and let $f:X\rightarrow Y$ be an embedding of $X$ into $Y$. We say that the embedding $f$ \emph{admits a retraction} if there is a retraction from $Y$ onto $f(X)$.
\end{definition}
In \cite{charatonik2}, W.~J.~Charatonik and J.~J.~Charatonik established two foundational results about embeddings of a Lelek fan or a Cantor fan into a Cantor fan, here we list them as Theorem \ref{C&C1} and Theorem \ref{C&C2}. 
\begin{theorem}\label{C&C1}
	Let $C$ be a Cantor fan and let $L$ be a Lelek fan. Then for any embedding $f:L\rightarrow C$, the embedding $f$ does not admit a retraction.
\end{theorem}
\begin{proof}
	The proof follows from \cite[Theorem 6, page 165]{charatonik2}.
\end{proof}
\begin{theorem}\label{C&C2}
	Let $C$ be a Cantor fan. Then for any embedding $f:C\rightarrow C$, the embedding $f$ admits a retraction.
\end{theorem}
\begin{proof}
	The proof follows from \cite[Theorem 6, page 165]{charatonik2}.
\end{proof}
Since this paper is about embeddings of Cantor fans or Lelek fans into Lelek fans or Cantor fans, we conclude this section by stating and proving Theorems \ref{gost} and \ref{sost} about topological properties of Cantor fans that are embedded into Lelek fans and vice-versa.

\begin{theorem}\label{gost}
Let $C$ be a Cantor fan and let $L$ be a Lelek fan. If $C\subseteq L$, then $C$ is nowhere dense in $L$. 	
\end{theorem}
\begin{proof}
Let $v$ be the top of $L$ (and of $C$). We show that $\Int_{L}(\Cl_{L}(C))=\emptyset$. Since $\Cl_{L}(C)=C$, it suffices to see that $\Int_{L}(C)=\emptyset$. Let $x_0\in C$ be any point and let $r_0>0$. We show that for the open ball $B_{L}(x_0,r_0)$ in $L$, centered in $x_0$ with radius $r_0$, 
$$
B_{L}(x_0,r_0)\not \subseteq C.
$$
 We consider the following possible cases.
\begin{enumerate}
	\item\label{motovilec} $x_0\not\in E(C)\cup\{v\}$. Let $r>0$ be such that $r<\min\{r_0, d(x_0,E(C)\cup \{v\})\}$. Since $x_0\in L$ and since $L$ is a Lelek fan, it follows that there is a sequence $(e_n)$ of endpoints of $L$ such that $\displaystyle \lim_{n\to\infty}e_n=x_0$. Let $n_0$ be a positive integer such that $e_{n_0}\in B_{L}(x_0,r)$. Then $e_{n_0}\in B_{L}(x_0,r_0)\setminus C$ and it follows that  $B_{L}(x_0,r_0)\not \subseteq C$.
	\item $x_0\in E(C)\cup\{v\}$. Let $x\in (B_{L}(x_0,r_0)\cap C)\setminus (E(C)\cup\{v\})$ and let $r>0$ be such that $B_{L}(x,r)\subseteq B_{L}(x_0,r_0)$. The argument from \ref{motovilec} shows that $B_{L}(x,r)\not \subseteq C$. It follows that $B_{L}(x_0,r_0)\not \subseteq C$. 
\end{enumerate}  
This completes the proof.
\end{proof}

\begin{theorem}\label{sost}
Let $C$ be a Cantor fan and let $L$ be a Lelek fan. If $L\subseteq C$, then $L$ is nowhere dense in $C$. 	
\end{theorem}
\begin{proof}
Let $v$ be the top of $C$ (and of $L$). We show that $\Int_{C}(\Cl_{C}(L))=\emptyset$. Since $\Cl_{C}(L)=L$, it suffices to see that $\Int_{C}(L)=\emptyset$. Let $x_0\in L$ be any point and let $r_0>0$. We show that 
$$
B_{C}(x_0,r_0)\not \subseteq L.
$$
Let $x\in (L\cap B_{C}(x_0,r_0))\setminus E(C)$ and let $(e_n)$ be a sequence in $E(L)$ such that $\displaystyle \lim_{n\to\infty}e_n=x$. Note that there is a positive integer $n_0$ such that for each $n\geq n_0$, $e_n\not\in E(C)$. For each positive integer $n\geq n_0$, let 
$$
x_n\in B_{C}\left(e_n,\frac{1}{2^n}\right)\setminus L
$$
It follows that $(x_n)$ is a sequence in $C\setminus L$ such that $\displaystyle \lim_{n\to\infty}x_n=x$. Therefore, there is a positive integer $m$ such that $x_m\in B_{C}(x_0,r_0)$. Since $x_m\not \in L$, it follows that $B_{C}(x_0,r_0)\not \subseteq L$.
\end{proof} 
\begin{observation}
	Let $L$ be a Lelek fan and let $C$ be a Cantor fan. Note that:
	\begin{enumerate}
		\item There are embeddings $f$ of $L$ into $L$ such that $f(L)$ is nowhere dense in $L$ and there are embeddings $f$ of $L$ into $L$ such that $f(L)$ is not nowhere dense in $L$. If the embedding $f$ is a homeomorphism, then $f(L)=L$, and, this is the only case when $f(L)$ is dense in $L$.
		\item There are embeddings $f$ of $C$ into $C$ such that $f(C)$ is nowhere dense in $C$ and there are embeddings $f$ of $C$ into $C$ such that $f(C)$ is not nowhere dense in $C$. If the embedding $f$ is a homeomorphism, then $f(C)=C$, and, this is the only case when $f(C)$ is dense in $C$. 
	\end{enumerate} 
\end{observation}
\section{Embeddings of the Cantor fan into the Lelek fan and retractions}\label{s2}
In this section, we prove that every embedding of the Cantor fan into the Lelek fan can be accompanied by a retraction from the Lelek fan onto the image of the Cantor fan.
\begin{theorem}\label{I1}
	Let $C$ be a Cantor fan and let $L$ be a Lelek fan. Then for any embedding $f:C\rightarrow L$, the embedding $f$ admits a retraction.
\end{theorem}
\begin{proof}
	Let $f:C\rightarrow L$ be an embedding of the Cantor fan $C$ into the Lelek fan $L$. To show that there is a retraction from $L$ onto $f(C)$, let $e:L\rightarrow C$ be an embedding of the Lelek fan $L$ into the Cantor fan $C$. Note that $e(f(C))$ is a Cantor fan in the Cantor fan $C$. By Theorem \ref{C&C2}, there is  a retraction from $C$ onto $e(f(C))$. So, let $p:C\rightarrow e(f(C))$ be a retraction from $C$ onto $e(f(C))$. Then we define the function $r:L\rightarrow f(C)$ as follows. For each $x\in L$, let 
	$$
	r(x)=\varphi(p(e(x))),
	$$
	where $\varphi:e(f(C))\rightarrow f(C)$ is defined by $\varphi(t)\in e^{-1}(t)$ for each $t\in e(f(C))$. Note that the preimage $e^{-1}(t)$ consists of a single element for each $t\in e(f(C))$. Therefore, since $e$ is an embedding, $\varphi$ is a well-defined continuous function. Since $r$ is a composition of three continuous functions, it is itself a continuous function. To prove that $r$ is a retraction from $L$ onto $f(C)$, let $x\in f(C)$. Then $e(x)\in e(f(C))$ and, since $p$ is a retraction from $C$ onto $e(f(C))$, it follows that $p(e(x))=e(x)$. Also, note that $\varphi(e(x))\in e^{-1}(e(x))$ and that $e^{-1}(e(x))=\{x\}$. Therefore,
	$$
	r(x)=\varphi(p(e(x)))=\varphi(e(x))=x.
	$$
	This completes the proof.
\end{proof}

\section{Embeddings of the Lelek fan into the Lelek fan and retractions}\label{s3}
Let $X$ and $Y$ be Lelek fans. In this section, we show that there exist embeddings \( f: X \to Y \) that admit retractions, as well as embeddings that do not. We further investigate this case by identifying specific conditions on the embedding \( f \) under which a retraction from \( Y \) onto \( f(X) \) does exist. We start with the following easy observation.

\begin{observation}\label{trama}
	Let $L$ be a Lelek fan and let $f:L\rightarrow L$ be a homeomorphism (note that in this case, $f(L)=L$). Then $f$ is an embedding of $L$ into $L$. Also, let $r:L\rightarrow f(L)$ be defined by $r(x)=x$ for any $x\in L$. Then $r$ is a retraction from $L$ onto $f(L)$. 
\end{observation}
Let \( L \) denote the Lelek fan. As seen in Observation \ref{trama}, any homeomorphism on \( L \) is an embedding of \( L \) into itself that admits a retraction. We call such an embedding a trivial embedding of \( L \) into itself. This naturally leads to the following question: are there any other, nontrivial embeddings of \( L \) into itself that also admit retractions? Among other things, Theorem \ref{I2} provides an answer to this question. 
\begin{theorem}\label{I2}
	Let $L$ be a Lelek fan. Then 
	\begin{enumerate}
		\item\label{1} there is an embedding $f:L\rightarrow L$ (other than a homeomorphism) of the Lelek fan $L$ into the Lelek fan $L$ such that $f$ admits a retraction.
		\item\label{2} there is an embedding $f:L\rightarrow L$ of the Lelek fan $L$ into the Lelek fan $L$ such that $f$ does not admit a retraction.
			\end{enumerate} 
\end{theorem}
\begin{proof}
	To prove \ref{1}, let $L_1$ and $L_2$ be lelek fans in the plane $\mathbb R^2$ such that
	\begin{enumerate}
		\item $L_1\cap L_2=\{(0,0)\}$ and 
		\item $(0,0)$ is the top of both fans, $L_1$ and $L_2$.
	\end{enumerate}
	Note that $L_1\cup L_2$ is again a Lelek fan, therefore, it is homeomorphic to $L$. Let $h:L_1\cup L_2\rightarrow L$ be a homeomorphism. Also, let $\varphi:L\rightarrow h(L_1)$ be a homeomorphism and let $f:L\rightarrow L$ be defined by $f(x)=\varphi(x)$ for any $x\in L$. Then $f$ is an embedding of the Lelek fan $L$ into the Lelek fan $L$. Next, let $r:L\rightarrow f(L)$ be defined by 
	$$
	r(x)=\begin{cases}
				x\text{;} & x\in h(L_1) \\
				\varphi(0,0)\text{;} & x\in h(L_2)
			\end{cases}
	$$
	for any $x\in L$. Note that $r$ is a continuous function from $L$ to $f(L)$, and, since $f(L)=\varphi(L_1)$, it follows that for each $x\in f(L)$, $r(x)=x$. Therefore, $r$ is a retraction from $L$ to $f(L)$. 
	
	To prove \ref{2}, let $C$ be a Cantor fan and let $e:C\rightarrow L$ and $g:L\rightarrow C$ be embeddings. Also, let $f:L\rightarrow L$ be defined by $f(x)=e(g(x))$ for any $x\in L$. Then $f$ is an embedding of the Lelek fan $L$ into the Lelek fan $L$. Suppose that there is a retraction from $L$ onto $f(L)$. Let $r:L\rightarrow f(L)$ be a retraction. Note that $e(g(L))$ is a Lelek fan, embedded into the Cantor fan $e(C)$. Let $p:e(C)\rightarrow e(g(L))$ be defined by 
	$$
	p(x)=r(x)
	$$
	for any $x\in e(C)$. Since $p$ is the restriction of the retraction $r$ to $e(C)$, it is a continuous function. Since for each $x\in e(g(L))$, $p(x)=r(x)=x$, it follows that $p$ is a retraction from the Cantor fan $e(C)$ onto the Lelek fan $e(g(L))$. This contradicts Theorem \ref{C&C1}. Therefore, in this case, there are no retractions from $L$ onto $f(L)$.
\end{proof}

The following is a corollary of the second part of the above proof.
\begin{theorem}\label{jjuurree}
	Let $L$ be a Lelek fan, let $C$ be a Cantor fan and let $f:L\rightarrow L$ be an embedding such that $f(L)\subseteq C\subseteq L$. Then $f$ does not admit a retraction. 
\end{theorem}
\begin{proof}
	Theorem follows directly from the second part of the proof of Theorem \ref{I2}.
\end{proof}

 \subsection{Wedges in Lelek fans}
\begin{figure}[h!]
	\centering
		\includegraphics[width=20em]{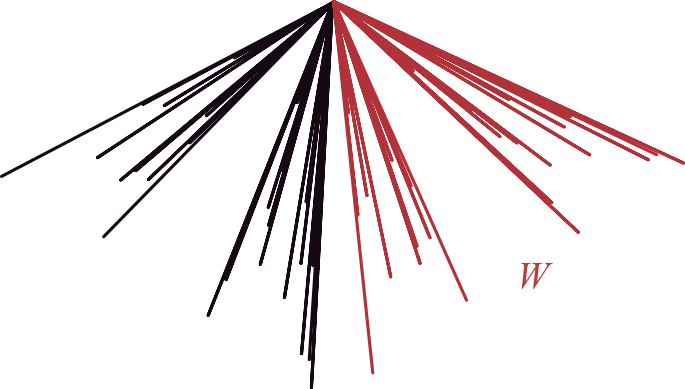}
	\caption{A wedge $W$ in a Lelek fan}
	\label{Fig2}
\end{figure}

Next, we present a  class of embeddings of the Lelek fan into itself that admit retractions; see Theorem \ref{jura}. The images of these embeddings (copies of the Lelek fan within $L$)  will be referred to as wedges in $L$. This class is formally described in Definition \ref{likja}. 
 \begin{definition}\label{likja}
 	Let $L$ be a Lelek fan with top $v$ and let  $W$ be a subcontinuum of $L$. If $W$ and $(L\setminus W)\cup \{v\}$ are both Lelek fans, then we say that $W$ is \emph{a wedge in $L$}; See Figure \ref{Fig2}, where a wedge $W$ in a Lelek fan is presented.
 \end{definition}
 \begin{observation}
	Let $L$ be a Lelek fan with top $v$ and let $W$ be a wedge in $L$. Then $v\in W$. 
\end{observation}
\begin{theorem}\label{jura}
Let $L$ be a Lelek fan with top $v$ and let $W$ be a wedge in $L$. Then there is a retraction from $L$ to $W$. 	
\end{theorem}
\begin{proof}
	Let $r:L\rightarrow W$ be defined by 
	$$
	r(x)=\begin{cases}
				x\text{;} & x\in W \\
				v\text{;} & x\in (L\setminus W)\cup \{v\}
			\end{cases}
	$$
	for any $x\in L$. Since $W$ and $(L\setminus W)\cup \{v\}$ are both closed in $L$ and since $$W\cap ((L\setminus W)\cup \{v\})=\{v\},$$ it follows that $r$ is a well-defined continuous function. It follows that $r$ is a retraction from $L$ to $W$.
\end{proof}
The following theorem is a characterization of a wedge in a Lelek fan. 
\begin{theorem}
	Let $L$ be a Lelek fan with top $v$ and let $W$ be a subset of $L$. The following statements are equivalent:
	\begin{enumerate}
		\item\label{eca} $W$ is a wedge in $L$.
		\item\label{ace} $W$ has the following properties:
	\begin{enumerate}
		\item $W\neq L$,
		\item $W$ is closed in $L$, 
		\item $W\setminus \{v\}$ is open in $L$, and
		\item for each leg $A\in \mathcal L(L)$,
		$$
		A\cap (W\setminus \{v\})\neq \emptyset ~~~ \Longrightarrow ~~~ A\subseteq W.
		$$
		\end{enumerate} 
		\end{enumerate} 
	\end{theorem}
\begin{proof}
	The implication from \ref{eca} to \ref{ace} follows from the definition of a wedge. To prove the implication from \ref{ace} to \ref{eca}, suppose that \begin{enumerate}
		\item[(a)] $W\neq L$,
		\item[(b)] $W$ is closed in $L$, 
		\item[(c)] $W\setminus \{v\}$ is open in $L$, and
		\item[(d)] for each leg $A\in \mathcal L(L)$,
		$$
		A\cap (W\setminus \{v\})\neq \emptyset ~~~ \Longrightarrow ~~~ A\subseteq W.
		$$
		\end{enumerate}
		We show that $W$ is a wedge in $L$ in the following steps.
		\begin{itemize}
			\item Step 1. Let $V=(L\setminus W)\cup \{v\}$. We show that 
		\begin{enumerate}
		\item[(i)] $V\neq L$,
		\item[(ii)] $V$ is closed in $L$, 
		\item[(iii)] $V\setminus \{v\}$ is open in $L$, and
		\item[(iv)] for each leg $A\in \mathcal L(L)$,
		$$
		A\cap (V\setminus \{v\})\neq \emptyset ~~~ \Longrightarrow ~~~ A\subseteq V.
		$$
		\end{enumerate} 		
	Since $W\neq L$, it follows that $V\neq L$. Since $W\setminus \{v\}$ is open in $L$, it follows that $L\setminus (W\setminus \{v\})$ is closed in $L$. Note that
	$$
	L\setminus (W\setminus \{v\})=(L\setminus W)\cup \{v\}=V. 
	$$
	Therefore, $V$ is closed in $L$. Also, note that $V\setminus \{v\}$ is open in $L$. Finally, let $A\in \mathcal L(L)$ be such a leg that $A\cap (V\setminus \{v\})\neq \emptyset$. To see that $A\subseteq V$, let $x\in A$. If $x=v$, then $x\in V$. Suppose next that $x\neq v$. If $x\not \in V$, then $x\in A\cap (W\setminus \{v\})$ and, since $W$ is a wedge in $L$, it follows that $A\subseteq W$, which is a contradiction. Therefore, $x\in V$ and $A\subseteq V$ follows. 
	\item Step 2. Let $\mathcal A=\{A\in \mathcal L(L) \ | \ A\subseteq W\}$. We show that $W=\bigcup_{A\in \mathcal A}A$. Note that $W\supseteq \bigcup_{A\in \mathcal A}A$ since for each $A\in \mathcal A$, $A\subseteq W$. To see that $W=\bigcup_{A\in \mathcal A}A$, let $x\in W$. First, suppose that $x=v$. Then $x\in \bigcup_{A\in \mathcal A}A$ since for each $A\in \mathcal A$, $x\in A$. Next, assume that $x\neq v$. Then $x\in A\setminus \{v\}$ for some leg $A\in \mathcal L(L)$. Let $A$ be such a leg. Since $W$ is a wedge in $L$, it follows that $A\subseteq W$. Therefore, $x\in \bigcup_{A\in \mathcal A}A$. This proves that $W=\bigcup_{A\in \mathcal A}A$. 
\item Step 3. We show that $W$ is a continuum. Let $\mathcal A=\{A\in \mathcal L(L) \ | \ A\subseteq W\}$. By Step 2, $W=\bigcup_{A\in \mathcal A}A$. Since $W$ is a union of arcs in $L$, each of them containing the point $v$, it follows that $W$ is connected. It follows from the definition of a wedge in $L$ that $W$ is compact, therefore, $W$ is a continuum. 
\item Step 4.   We prove that $W$ is a Lelek fan. By Step 3, $W$ is a continuum. Therefore, since it is a subcontinuum of the smooth fan $L$, it is itself a smooth fan. Let $\mathcal A=\{A\in \mathcal L(L) \ | \ A\subseteq W\}$. By Step 2, $W=\bigcup_{A\in \mathcal A}A$. To see that $W$ is a Lelek fan, we prove that $\Cl (E(W))=W$. Note that $\Cl (E(W))\subseteq W$. To see that $\Cl (E(W))\supseteq W$, let $x\in W$. Without loss of generality, suppose that $x\neq v$. We prove that there is a sequence $(e_n)$ in $E(W)$ such that $\displaystyle\lim_{n\to \infty}e_n=x$. Suppose that for each sequence $(e_n)$ in  $E(W)$, $\displaystyle\lim_{n\to \infty}e_n\neq x$. Since $x\in L$ and since $L$ is a Lelek fan, there is a sequence $(e_n)$ in $E(L)$ such that $\displaystyle\lim_{n\to \infty}e_n=x$. Let $(e_n)$ be such a sequence. Since for each sequence $(e_n)$ in  $E(W)$, $\displaystyle\lim_{n\to \infty}e_n\neq x$, it follows that for each subsequence $(e_{i_n})$ of the sequence $(e_n)$, $(e_{i_n})$ is not a sequence in $W$. Therefore, $(L\setminus W)\cup \{v\}$ is not closed in $L$. This contradicts the properties of $(L\setminus W)\cup \{v\}$ from Step 1. 
\item Step 5. A similar argument as in Step 4 shows that also $(L\setminus W)\cup \{v\}$ is a Lelek fan (note that it follows from Step 1 of this proof that $W$ and $(L\setminus W)\cup \{v\}$ have similar properties that lead to the same conclusion that $(L\setminus W)\cup \{v\}$ is a Lelek fan).
	\end{itemize}
	This proves that $W$ is a wedge in $L$. 
\end{proof}

\subsection{Cuts of Lelek fans}
In this section, we introduce the concept of cuts in Lelek fans and show that these cuts allow embeddings of a Lelek fan into itself that admit retractions from the fan onto the image of the embedding. We begin with the following lemma.
\begin{lemma}\label{abrakadabra}
	Let $X$ be a metric space, let $Y$ and $Z$ be subspaces of $X$ such that $Z\cap Y\neq \emptyset$ and $Z\cap (X\setminus Y)\neq \emptyset$. If $Y$ is a retract of $X$, then $Y\cap Z$ is a retract of $Z$.  
\end{lemma}
\begin{proof}
	Suppose that $Y$ is a retract of $X$. We prove that $Y\cap Z$ is a retract of $Z$. Let $r:X\rightarrow Y$ be a retraction and let $i:Z\rightarrow X$ be the inclusion function; i.e., for each $z\in Z$, $i(z)=z$. Also, let $R:Z\rightarrow Y\cap Z$ be defined by 
	$$
	R=r\circ i.
	$$
	Since $R$ is composed of continuous functions $r$ and $i$, it is itself a continuous function. Let $y\in Y\cap Z$. Then 
	$$
	R(y)=(r\circ i)(y)=r(i(y))=r(y)=y,
	$$
	and, therefore, $R$ is a retraction from $Z$ to $Y\cap Z$. 
\end{proof}
\begin{definition}
	Let $X$ be a metric space. Then we define the family $\mathcal R(X)$ by
	$$
	\mathcal R(X)=\{Y\subseteq X \ | \ \textup{there is a retraction from } X \textup{ to } Y\}.
	$$
	Also, let $Z\subseteq X$. Then we define $\mathcal C(Z,X)$ by
	$$
	\mathcal C(Z,X)=\{Y\in \mathcal R(X) \ | \ Z\cap Y\neq \emptyset \textup{ and } Z\cap (X\setminus Y)\neq \emptyset\}.
	$$
	For each $Y\in \mathcal C(Z,X)$ we say that $Y$ is \emph{a cut of $Z$ in $X$}.
\end{definition}
\begin{theorem}\label{cuts}
	Let $X$ be a metric space and let $Y$ and $Z$ be subspaces of $X$. If $Y$ is a cut of $Z$ in $X$, then  $Y\cap Z$ is a retract of $Z$.
\end{theorem}
\begin{proof}
	Let $Y$ be a cut of $Z$ in $X$. Then $Z\cap Y\neq \emptyset$ and $Z\cap (X\setminus Y)\neq \emptyset$. Also, it follows from $Y\in \mathcal R(X)$ that $Y$ is a retract of $X$. By Lemma \ref{abrakadabra}, $Y\cap Z$ is a retract of $Z$. 
\end{proof}
Note that Lelek fans are planar continua. The following corollary is about cuts of Lelek fans in $\mathbb R^2$.
\begin{corollary}\label{juzina}
	Let $L$ be a Lelek fan in $\mathbb R^2$ and let $f:L\rightarrow L$ be an embedding. If there is a cut $Y$ of $L$ in $\mathbb R^2$ such that $Y\cap L=f(L)$, then $f$ admits a retraction. 
\end{corollary}
\begin{proof}
	The corollary follows directly from Theorem \ref{cuts}.
\end{proof}
The following example demonstrates how Corollary \ref{juzina} may be used.
\begin{example}\label{pans}
	Let $X=\mathbb R^2$ and let $F$ be the Cantor fan $F=\bigcup_{c\in C}A_c$, where $C\subseteq [0,1]$ is the Cantor middle third set and for each $c\in C$, $A_c$ is the  convex segment in the plane from $(\frac{1}{2},0)$ to $(c,1)$. Also, let $L$ be a Lelek fan in $F$; see Figure \ref{Fig3} and let 
	$$
	\varepsilon \in \left(0,\max\left\{\frac{d((\frac{1}{2},0),(x,y))}{2} \ \Big | \ (x,y)\in L\right\}\right).
	$$  
	\begin{figure}[h!]
	\centering
		\includegraphics[width=14em]{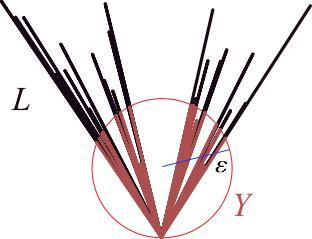}
	\caption{The Lelek fan $L$ from Example \ref{pans}}
	\label{Fig3}
\end{figure}
Let $Y=\Cl\Big(B\big((\frac{1}{2},\varepsilon),\varepsilon\big)\Big)$ and note that $Y$ is a cut of $L$ in $\mathbb R^2$. Also, note that $Y\cap L$ is a Lelek fan in $L$. Therefore, for any embedding $f:L\rightarrow L$ such that $f(L)=Y\cap L$, the embedding $f$ admits a retraction.   
\end{example}
\subsection{Combining different variants}
Here, we merge the concepts of wedges, cuts, and nowhere dense Lelek fans to derive additional results on self-embeddings of Lelek fans that either admit or fail to admit retractions. Theorem \ref{wer} illustrates how the notion of a wedge can be combined with that of a cut while Theorem \ref{wer1} deals with a situation where a nowhere dense Lelek fan is captured in a wedge.
\begin{theorem}\label{wer} 
	Let $L$ be a Lelek fan in $\mathbb R^2$ with top $v$, let $W_1$ and $W_2$ be wedges in $L$ such that $W_1\cup W_2=L$ and $W_1\cap W_2=\{v\}$, and let $Y\subseteq \mathbb R^2$ be such that 
	$Y\cap W_2$ is a Lelek fan. Then the following hold.
\begin{enumerate}
	\item $W_1\cup (Y\cap W_2)$ is a Lelek fan in $L$.
	\item If $Y$ is a cut of $W_2$ in $\mathbb R^2$, then $W_1\cup (Y\cap W_2)$ is a retract of $L$.
	\item $Y\cap W_2$ is a retract of $L$.
\end{enumerate}	
		\end{theorem}
\begin{proof}
	\begin{figure}[h!]
	\centering
		\includegraphics[width=17em]{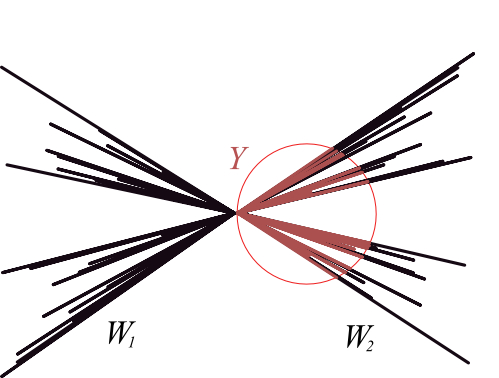}
	\caption{The Lelek fan $L$ from Theorem \ref{wer}}
	\label{Fig4}
\end{figure}
Since $W_1\cup (Y\cap W_2)$ is a union of two Lelek fans that intersect exactly in their tops, it is a Lelek fan in $L$. Let $r:W_2\rightarrow Y\cap W_2$ be a retraction and let $R:L\rightarrow W_1\cup (Y\cap W_2)$ be defined by 
$$
R(x)=\begin{cases}
				x\text{;} & x\in W_1 \\
				r(x)\text{;} & x\in W_2
			\end{cases}
$$
for each $x\in L$. Then $R$ is a retraction from $L$ onto $W_1\cup (Y\cap W_2)$. Additionally, the function $\rho:L\rightarrow Y\cap W_2$, defined by 
$$
\rho(x)=\begin{cases}
				v\text{;} & x\in W_1 \\
				R(x)\text{;} & x\in W_2
			\end{cases}
$$
for each $x\in L$, is a retraction.
\end{proof}

\begin{theorem}\label{wer1}
	Let $L$ be a Lelek fan with top $v$, let $W$ be a wedge in $L$, and let $L'$ be a Lelek fan in $L$ such that 
	\begin{enumerate}
	\item $L'\cap W$ is a Lelek fan and 
	\item for some Cantor fan $C$, $L'\cap W\subseteq C\subseteq W$.
\end{enumerate}
Then there there are no retractions from $L$ onto $L'$.	
		\end{theorem}
\begin{proof}
Suppose that there is a retraction $r$ from $L$ onto $L'$. Let $i:W\rightarrow L$ be the inclusion function; i.e., for each $x\in W$, $i(x)=x$. Also, let $R:W\rightarrow L'\cap W$ be defined by $R=r\circ i$.  Since $R$ is composed of continuous functions $r$ and $i$, it is itself a continuous function. Let $x\in L'\cap W$.  Then 
	$$
	R(x)=(r\circ i)(x)=r(i(x))=r(x)=x,
	$$
	and, therefore, $R$ is a retraction from $W$ to $L'\cap W$. This contradicts Theorem \ref{jjuurree}. Therefore,  there there are no retractions from $L$ onto $L'$.	
\end{proof}
\subsection{Confluent retractions}
In this section, we study continua that arise as confluent images of the Lelek fan.
It was shown in \cite{charatonik} that any such image is homeomorphic to the Lelek fan, provided it is not a one-point continuum.
This result plays a central role in the proof of our main theorem in this section, Theorem \ref{jeza}, where confluent retractions on the Lelek fan are studied.
We begin with the following definitions.
\begin{definition}
 Let $X$ and $Y$ be continua. A continuous  mapping $f: X \rightarrow  Y$ is said to be \emph{confluent}, if for each subcontinuum $Q$ of $f(X)$ and for each component $K$ of $f^{-1}(Q)$, $f(K)=Q$.
\end{definition}
\begin{observation}
	Note that the retraction $r$ used in the proof of Theorem \ref{jura} is confluent.
	\end{observation}

\begin{theorem}\label{jeza}
	Let $L$ be a Lelek fan and let $Y$ be a non-degenerate continuum. The following statements are equivalent.
	\begin{enumerate}
	    \item\label{tre} There is an embedding $e:Y\rightarrow L$ of $Y$ into $L$ that admits a confluent retraction from $L$ onto $e(Y)$.
		\item\label{uno} There is a confluent surjective mapping $f:L\rightarrow Y$.
		\item\label{nove} $Y$ is a Lelek fan. 
		\end{enumerate}
\end{theorem}
\begin{proof}
To prove the implication from \ref{tre} to \ref{uno}, let $e:Y\rightarrow L$ be an embedding of $Y$ into $L$,  let $r:L\rightarrow e(Y)$ be a confluent retraction, let $\varphi:e(Y)\rightarrow Y$ be a homeomorphism, and let $f=\varphi\circ r$. Then $f:L\rightarrow Y$ and it is a confluent surjection.  The proof of the implication from \ref{uno} to \ref{nove} can be found in \cite{charatonik}; it follows from \cite[Theorem 13]{charatonik}.  Next, we prove the implication  from \ref{nove} to \ref{tre}. Let $Y$ be a Lelek fan. Let $e:Y\rightarrow L$ be a homeomorphism. Also, let $r:L\rightarrow e(Y)$ be defined by 
	$$
	r(x)=x
	$$
	for any $x\in L$. Then $r$ is a retraction from $L$ to $e(Y)$. It follows from the definition of the retraction $r$ that it is a confluent mapping. Therefore, \ref{tre} follows. 
\end{proof}

\subsection{Monotone retractions}
In this section, we study continua that arise as monotone images of the Lelek fan.
It is a well-known fact that any such image is homeomorphic to the Lelek fan, if it is not a one-point continuum.
This result plays a central role in the proof of our main theorem in this section, Theorem \ref{jeza1}, where monotone retractions on the Lelek fan are studied.
We begin with the following definitions.  
\begin{definition}
 Let $X$ and $Y$ be continua. A continuous  mapping $f: X \rightarrow  Y$ is said to be \emph{monotone} if the preimage $f^{-1}(y)$ is connected for each point $y\in Y$.
\end{definition}
We  use Lemma \ref{Judy} in the proof of Theorem \ref{jeza1}.
\begin{lemma}\label{Judy}
	Let $X$ and $Y$ be continua and let $f:X\rightarrow Y$ be a monotone mapping. Then $f$ is confluent.
\end{lemma}
\begin{proof}
	See \cite[page 164]{charatonik}.
\end{proof}
\begin{theorem}\label{jeza1}
	Let $L$ be a Lelek fan and let $Y$ be a non-degenerate continuum. The following statements are equivalent.
	\begin{enumerate}
	\item\label{tre1} There is an embedding $e:Y\rightarrow L$ of $Y$ into $L$ that admits a monotone retraction from $L$ onto $e(Y)$.
		\item\label{uno1} There is a monotone surjective mapping $f:L\rightarrow Y$.
		\item\label{nove1} $Y$ is a Lelek fan. 
	\end{enumerate}
\end{theorem}
\begin{proof}
To prove the implication from \ref{tre1} to \ref{uno1}, let $e:Y\rightarrow L$ be an embedding of $Y$ into $L$,  let $r:L\rightarrow e(Y)$ be a monotone retraction, let $\varphi:e(Y)\rightarrow Y$ be a homeomorphism, and let $f=\varphi\circ r$. Then $f:L\rightarrow Y$ and it is a monotone surjection.  To prove the implication \ref{uno1} to \ref{nove1}, let $f:L\rightarrow Y$ be a monotone surjection. By Lemma \ref{Judy}, $f$ is a confluent surjection and it follows from Theorem \ref{jeza} that $Y$ is a Lelek fan.  Next, we prove the implication  from \ref{nove1} to \ref{tre1}. Let $Y$ be a Lelek fan. Let $e:Y\rightarrow L$ be a homeomorphism. Also, let $r:L\rightarrow e(Y)$ be defined by 
	$$
	r(x)=x
	$$
	for any $x\in L$. Then $r$ is a retraction from $L$ to $e(Y)$. It follows from the definition of the retraction $r$ that it is a monotone mapping. Therefore,  \ref{tre1} follows. 
\end{proof}
\subsection{Open retractions}
In this section, we give results about continua that arise as open images of the Lelek fan.
Since any open mapping is confluent, it follows that any such image is homeomorphic to the Lelek fan, provided it is not a one-point continuum.
This result plays a central role in the proof of our main theorem in this section, Theorem \ref{jeza11}, where results about open retractions on the Lelek fan are presented.
\begin{definition}
	Let $X$ and $Y$ be continua. A continuous  mapping $f: X \rightarrow  Y$ is said to be \emph{open}, if the image of any open subset of $X$ is an open subset of $f(X)$.
\end{definition}
\begin{observation}
	Note that the retraction $r$ used in the proof of Theorem \ref{jura} is monotone and it is not open. 
\end{observation}
We use the following lemma in the proof of Theorem \ref{jeza11}. 
\begin{lemma}\label{Judy2}
	Let $X$ and $Y$ be continua and let $f:X\rightarrow Y$ be an open mapping. Then $f$ is confluent.
\end{lemma}
\begin{proof}
	See \cite[Theorem 7.5]{W}.
\end{proof}
\begin{theorem}\label{jeza11}
	Let $L$ be a Lelek fan and let $Y$ be a non-degenerate continuum. The following statements are equivalent.
	\begin{enumerate}
	\item\label{tre11} There is an embedding $e:Y\rightarrow L$ of $Y$ into $L$ that admits an open retraction from $L$ onto $e(Y)$.
			\item\label{uno11} There is an open surjective mapping $f:L\rightarrow Y$.
		\item\label{nove11} $Y$ is a Lelek fan. 
	\end{enumerate}
\end{theorem}
\begin{proof}
To prove the implication from \ref{tre11} to \ref{uno11}, let $e:Y\rightarrow L$ be an embedding of $Y$ into $L$,  let $r:L\rightarrow e(Y)$ be an open retraction, let $\varphi:e(Y)\rightarrow Y$ be a homeomorphism, and let $f=\varphi\circ r$. Then $f:L\rightarrow Y$ and it is an open surjection.  To prove the implication \ref{uno11} to \ref{nove11}, let $f:L\rightarrow Y$ be an open surjection. By Lemma \ref{Judy}, $f$ is a confluent surjection and it follows from Theorem \ref{jeza} that $Y$ is a Lelek fan.  Next, we prove the implications  from \ref{nove11} to \ref{tre11}. Let $Y$ be a Lelek fan. Let $e:Y\rightarrow L$ be a homeomorphism. Also, let $r:L\rightarrow e(Y)$ be defined by 
	$$
	r(x)=x
	$$
	for any $x\in L$. Then $r$ is a retraction from $L$ to $e(Y)$. It follows from the definition of the retraction $r$ that it is an open  mapping. Therefore,  \ref{tre11} follows. 
	\end{proof}

	\section{Simple and semi-simple retractions}\label{s4}
	
	Let $L$ and $L'$ be Lelek fans such that $L'\subseteq L$. In this section we define a simple and a semi-simple retraction from $L$ onto $L'$ and study retraction $r:L\rightarrow L'$ that  can be factorized through a retraction from the fence of $L$ onto the fence of $L'$. We start with the following definition.

\begin{definition}
We always use $Y$ to denote a Cantor set. Also, we use $F$ to denote $F=Y\times [0,1]$ and $\sim$ to denote the equivalence relation on $F$, defined by 
$$
(x,t)\sim (y,s) ~~~  \Longleftrightarrow  ~~~   (x,t) = (y,s) \textup{ or } s=t=0.
$$
for all $(x,t),(y,s)\in F$. Then we define 
 $$
 C=F/_{\sim}.
 $$
  \end{definition}
  \begin{observation}
  	Note that $C$ is a Cantor fan. 
  \end{observation}
  \begin{definition}
We use $q:F\rightarrow C$ to denote the quotient mapping defined by $q(x,t)=[(x,t)]$ for each $(x,t)\in F$, where $[(x,t)]$ is the equivalence class of $(x,t)$ for each $(x,t)\in F$. 
 \end{definition}
 \begin{observation}
  	Note that for each $(x,t)\in F$, $q(x,t)=\{(x,t)\}$ if and only if $t\neq 0$. Also note that for each $x\in Y$,  $q(x,0)=Y\times \{0\}$.
  \end{observation}
  \begin{definition}
For each continuum $K$ in $C$, we use $F_K$ and $H_K$ to denote 
 $$
 F_K=q^{-1}(K) ~~~ \textup{ and } ~~~ H_K=\Cl_F(q^{-1}(K\setminus \{v\})).
 $$
 We call $F_K$ and $H_K$ \emph{the fences of the continuum $K$}.
 \end{definition}
  \begin{observation}
 	Note that $F_C=F$. Also, note that for each continuum $K$ in $C$, $H_K\subseteq F_K$.
 \end{observation}
\subsection{Simple retractions}\label{cebula}
 In this section we define a simple retraction from $L$ onto $L'$ and show in Corollary \ref{miselovka} that for any two Lelek fans  $L$ and $L'$ such that $L'\subseteq L$, and for any retraction $r:L\rightarrow L'$, the following statements are equivalent:
\begin{enumerate}
	\item The retraction $r$ can be factorized through a retraction from the fence $F_L$ onto the fence $F_{L'}$. 
	\item The retraction $r$ is a simple retraction from $L$ onto $L'$. 
\end{enumerate}

 \begin{definition}
 	We use $J$ and $K$ to denote fans in $C$ such that $K\subseteq J$. We also use $v$ to denote the ramification point in $C$. 
 \end{definition}
 \begin{observation}
 	Note that $v$ is also the ramification point in $J$ as well as in $K$. 
 \end{observation}
 \begin{definition}
 	For each $A\in \mathcal L(J)$, we use $L_A$ to denote the non-degenerate component of $q^{-1}(A)$. We also use $\mathcal L(F_J)$ to denote the set 
 	$$
 	\mathcal L(F_J)=\{L_A \ | \ A\in \mathcal L(J)\}. 
 	$$
\end{definition}
  \begin{observation}
 	Note that for each $H\in \mathcal L(F_J)$, $q(H)\in \mathcal L(J)$.
 \end{observation}
  \begin{definition}
 	For each $A\in \mathcal L(J)$, let $v_A$ be defined by 
 	$$
 	v_A\in q^{-1}(v)\cap L_A. 
 	$$
 	 \end{definition}

 \begin{definition}
 	Suppose that there is a retraction from $J$ onto $K$ and let $r:J\rightarrow K$ be a retraction. Then we say that $r$ is \emph{a simple retraction} if for each $A\in \mathcal L(J)$, 
 	$$
 	r(A)\subseteq A. 
 	$$
 \end{definition}

 \begin{theorem}\label{FrEd}
 	Suppose that there is a retraction from $J$ onto $K$ and let $r:J\rightarrow K$ be a retraction. The following statements are equivalent. 
 	\begin{enumerate}
 	\item\label{sina2} The retraction $r$ is simple.
 	\item\label{sina1} There is a retraction $R:F_J\rightarrow F_{K}$ such that for each $x\in F_J$, 
 		$$
 		q(R(x)) = r(q(x)).
 		$$
 	\end{enumerate} 
 \end{theorem}
 \begin{proof}
 To prove the implication from \ref{sina2} to \ref{sina1}, let $r:J\rightarrow K$ be a simple retraction. Let $x\in F_J$. If there is $A_0\in \mathcal L(J)$ such that $x\in L_{A_0}$, then the set $L_{A_0}\cap q^{-1}(r(q(x)))$ consists of a single point. So, for each  $x\in F_J$ such that there is $A_0\in \mathcal L(J)$ such that $x\in L_{A_0}$, let $y_x\in L_{A_0}\cap q^{-1}(r(q(x)))$. We define $R(x)$  by 
$$
R(x)=\begin{cases}
				y_x\text{;} & \textup{there is } A_0\in \mathcal L(J) \textup{ such that } x\in L_{A_0}\\
				x\text{;} & \textup{otherwise}
			\end{cases}
$$
for each $x\in F_J$. Note that $R$ is a well-defined function.

 Next, we prove that for each $x\in F_{K}$, $R(x)=x$. Let $x\in F_{K}$. We consider the following cases.
 \begin{enumerate}
 	\item[(a)] There is $A_0\in \mathcal L(J)$ such that $x\in L_{A_0}$. Choose and fix such a leg $A_0$.  It follows that $q(x)\in A_0\cap K$. Therefore, $r(q(x))=q(x)$ and 
 	$$
 	L_{A_0}\cap q^{-1}(r(q(x)))=L_{A_0}\cap q^{-1}(q(x))=\begin{cases}
				\{x\}\text{;} & q(x)\neq v\\
				\{v_{A_0}\}\text{;} & q(x)=v.
			\end{cases}
 	$$
 	If $q(x)\neq v$, then $R(x)=y_x=x$. If $q(x)=v$, then $x=v_{A_0}$ and it follows that $R(x)=y_x=v_{A_0}=x$.
 	\item[(b)] For each $A\in \mathcal L(J)$, $x\not\in L_{A}$. Then $R(x)=x$ and we are done.
 \end{enumerate}
 Next, we prove that for each $x\in F_J$, $q(R(x)) = r(q(x))$. Let $x\in F_J$. We consider the following cases.
 \begin{enumerate}
 	\item[(a)] There is $A_0\in \mathcal L(J)$ such that $x\in L_{A_0}$. Choose and fix such a leg $A_0$.  It follows that $q(x)\in A_0$. Since $R(x) = y_x\in L_{A_0}\cap q^{-1}(r(q(x)))$, it follows that 
 	$$
 	q(R(x))=q(y_x)=r(q(x))
 	$$
 	and we are done.
 	\item[(b)] For each $A\in \mathcal L(J)$, $x\not\in L_{A}$. Then $q(x)=v$ and it follows that $q(R(x))=q(x)=v$ and $r(q(x))=r(v)=v$ and $q(R(x)) = r(q(x))$ follows.
 \end{enumerate}
Finally, we prove that $R$ is continuous. Let $x\in F_J$ and let $(x_n)$ be a sequence in $F_J$ such that $\displaystyle\lim_{n\to\infty}x_n=x$. We show that $\displaystyle \lim_{n\to\infty}R(x_n)=R(x)$. We consider the following cases.
\begin{enumerate}
	\item[(a)] There is $A_0\in \mathcal L(J)$ such that $x\in L_{A_0}$. Choose and fix such a leg $A_0$. We consider the following subcases.
	\begin{enumerate}
	\item[(i)] $x\neq v_{A_0}$. Note that in this case there is a positive integer $n_0$ such that for each positive integer $n$,
	$$
	n\geq n_0 ~~~ \Longrightarrow ~~~ x_n\not\in Y\times \{0\}.
	$$
	Without loss of generality suppose that $n_0=1$. Note that $q(R(x))=\{R(x)\}=r(q(x))$ and that for each positive integer $n$, $q(R(x_n))=\{R(x_n)\}=r(q(x_n))$. Since $r$ and $q$ are continuous, it follows that 
	$$
	\lim_{n\to\infty}r(q(x_n))=r(q(x)).
	$$
	Therefore,
	$$
	\lim_{n\to\infty}R(x_n)=R(x).
	$$
	\item[(ii)] $x=v_{A_0}$. Then $R(x)=x$. We consider the following cases.
	\begin{itemize}
	\item There is a positive integer $n_0$ such that for each positive integer $n$,
	$$
	n\geq n_0 ~~~ \Longrightarrow ~~~ \textup{there is } A_n\in \mathcal L(J) \textup{ such that } x_n\in L_{A_n}.
	$$
	Without loss of generality suppose that $n_0=1$. For each positive integer $n$, let $A_n\in \mathcal L(J)$ be such that $x_n\in L_{A_n}$. Note that $R(x)=v_{A_0}$ and that for each positive integer $n$, $d(R(x_n),v_{A_0})\leq d(x_n,v_{A_0})$. Therefore, since $\displaystyle \lim_{n\to\infty}x_n=x$ ($=v_{A_0}$) it follows that 
	$$
	\lim_{n\to\infty}R(x_n)=R(x).
	$$
	\item For each positive integer $n$, there is a positive integer $i_n>n$ such that for each $A\in \mathcal L(J)$, $x_{i_n}\not\in L_{A_{i_n}}$. 
		Without loss of generality suppose that for each positive integer $n$, $x_{n}\not\in L_{A_{n}}$. Then for each positive integer $n$, let $R(x_n)=x_n$. It follows that 
	$$
	\lim_{n\to\infty}R(x_n)=\lim_{n\to\infty}x_n=x=R(x).
	$$
	\end{itemize}
	\end{enumerate}
	\item[(b)] For each $A\in \mathcal L(J)$, $x\not\in L_{A}$. It follows that $R(x)=x$. We consider the following subcases.
	\begin{enumerate}
		\item[(i)] There is a positive integer $n_0$ such that for each positive integer $n$,
	$$
	n\geq n_0 ~~~ \Longrightarrow ~~~ \textup{there is } A_n\in \mathcal L(J) \textup{ such that } x_n\in L_{A_n}.
	$$
	Without loss of generality suppose that $n_0=1$. For each positive integer $n$, let $A_n\in \mathcal L(J)$ be such that $x_n\in L_{A_n}$. Note that $R(x)=x$ and that for each positive integer $n$, 
	$$
	d(R(x_n),R(x))=d(R(x_n),x)\leq d(x_n,x).
	$$
	 Therefore, since $\displaystyle \lim_{n\to\infty}x_n=x$  it follows that 
	$$
	\lim_{n\to\infty}R(x_n)=R(x).
	$$
	\item[(ii)] For each positive integer $n$, there is a positive integer $i_n>n$ such that for each $A\in \mathcal L(J)$, $x_{i_n}\not\in L_{A_{i_n}}$. 
		Without loss of generality suppose that for each positive integer $n$, $x_{n}\not\in L_{A_{n}}$. Then for each positive integer $n$, let $R(x_n)=x_n$. It follows that 
	$$
	\lim_{n\to\infty}R(x_n)=\lim_{n\to\infty}x_n=x=R(x).
	$$
	\end{enumerate}
\end{enumerate}
 This proves that $R$ is continuous.
 	
 		To prove the implication from \ref{sina1} to \ref{sina2}, let $R:F_J\rightarrow F_{K}$ be a retraction such that for each $x\in F_L$, 
 		$$
 		q(R(x)) = r(q(x))
 		$$
	and suppose that $r$ is not simple.  Then there are $A,B\in \mathcal L(J)$ such that $A\neq B$ and $r(A)\cap (B\setminus \{v\})\neq \emptyset$. Such legs $A$ and $B$ do exist since $r$ is not simple. Since $R$ is a retraction from $F_J$ to $F_{K}$ and since $Y\times \{0\}\subseteq F_{K}$, it follows that $R(L_A)\subseteq L_A$. Therefore,
	$$
	q(R(L_A))\subseteq q(L_A)=A.
	$$
	But
	$$
	q(R(L_A))=r(q(L_A))=r(A)\not \subseteq A
	$$
	since $r(A)\cap (B\setminus \{v\})\neq \emptyset$. This is a contradiction. Therefore, $r$ is a simple retraction.
 \end{proof}

 \begin{observation}
 	Note that Theorem \ref{FrEd} does not hold if $F_K$ and $F_J$ are replaced by $H_K$ and $H_J$. See Example \ref{juro} for details.
 \end{observation}
 \begin{corollary}\label{miselovka}
 	Let $L$ and $L'$ be Lelek fans in $F$ such that $L'\subseteq L$, suppose that there is a retraction from $L$ onto $L'$ and let $r:L\rightarrow L'$ be a retraction. The following statements are equivalent. 
 	\begin{enumerate}
 	\item\label{sina2} The retraction $r$ is simple.
 	\item\label{sina1} There is a retraction $R:F_L\rightarrow F_{L'}$ such that for each $x\in F_L$, 
 		$$
 		q(R(x)) = r(q(x)).
 		$$
 	\end{enumerate} 
 \end{corollary}
 \begin{proof}
 	Let $J=L$ and $K=L'$ and use Theorem \ref{FrEd}.
 \end{proof}
 
\subsection{Semi-simple retractions}
 In the final part of this paper, we explore a weakening of the notion of simple retractions -- we introduce the class of \emph{semi-simple retractions}. These retractions preserve certain leg structures in a relaxed sense, while still interacting naturally with corresponding retractions on the associated fences $H_L$ and $H_{L'}$. This leads us to new examples and counterexamples, and culminates in an open problem concerning the conditions under which such compatible retractions exist.

 \begin{definition}
 	We use $J$ and $K$ to denote fans in $C$ such that $K\subseteq J$. We also use $v$ to denote the ramification point in $C$. 
 \end{definition}
 \begin{definition}
 	Suppose that there is a retraction from $J$ onto $K$ and let $r:J\rightarrow K$ be a retraction. We say that $r$ is \emph{semi-simple} if for each $A\in \mathcal L(J)$ there is $B\in \mathcal L(J)$ such that 
 	$$
 	r(A)\subseteq B.
 	$$
 \end{definition}
 \begin{observation}
 	Note that every simple retraction is also a semi-simple retraction.
 \end{observation}

  \begin{definition}
 	For each $A\in \mathcal L(J)$, we use $L_A$ to denote the non-degenerate component of $q^{-1}(A)$. We also use $\mathcal L(H_J)$ to denote the set 
 	$$
 	\mathcal L(H_J)=\{L_A \ | \ A\in \mathcal L(J)\}. 
 	$$
\end{definition}
  \begin{observation}
 	Note that for each $H\in \mathcal L(H_J)$, $q(H)\in \mathcal L(J)$.
 \end{observation}
 \begin{theorem}\label{FrEdek}
 	Let $r:J\rightarrow K$ be a retraction and suppose that there is a retraction $R:H_J\rightarrow H_{K}$ such that for each $x\in H_J$, 
 		$$
 		q(R(x)) = r(q(x)).
 		$$
 		Then $r$ is semi-simple. 
 \end{theorem}
 \begin{proof}
 Suppose that $R:H_J\rightarrow H_{K}$ is such a retraction that for each $x\in H_J$, 
 		$$
 		q(R(x)) = r(q(x)).
 		$$
 Also, suppose that $r$ is not semi-simple.  Let $A_0,B_{0},C_{{0}}\in \mathcal L(J)$ be such that $B_0\neq C_0$ and
 	$$
 	r(A_0)\cap (B_0\setminus \{v\})\neq \emptyset ~~~ \textup{ and } ~~~ r(A_0)\cap (C_0\setminus \{v\})\neq \emptyset.
 	$$
 	Such legs $A_0$, $B_0$ and $C_0$ do exist since $r$ is not simple. Since $(L_{B_0}\cap H_{K})\cap R(L_{A_0})\neq \emptyset$ and $(L_{C_0}\cap H_{K})\cap R(L_{A_0})\neq \emptyset$, it follows that $R(L_{A_0})$ is not connected. Since $R$ is continuous and since $L_{A_0}$ is an arc, it follows that $R(L_{A_0})$ is connected. This is a contradiction.  This proves that the retraction $r$ is semi-simple.
 \end{proof}
\begin{corollary}
	Let $L$ and $L'$ be Lelek fans in $F$ such that $L'\subseteq L$, suppose that there is a retraction from $L$ onto $L'$ and let $r:L\rightarrow L'$ be a retraction. Also, suppose that there is a retraction $R:H_L\rightarrow H_{L'}$ such that for each $x\in H_L$, 
 		$$
 		q(R(x)) = r(q(x)).
 		$$
 		Then $r$ is semi-simple. 
\end{corollary}
\begin{proof}
	Let $J=L$ and $K=L'$ and use Theorem \ref{FrEdek}.
\end{proof}
 Next, we discuss the converse of Theorem \ref{FrEdek}.  First, we give the following observation.
 \begin{observation}
 	Note that if $r$ is simple, then a slight modification of the proof of Theorem \ref{FrEd} shows that there is a retraction $R:H_{J}\rightarrow H_{K}$ such that for each $x\in H_L$, 
 		$$
 		q(R(x)) = r(q(x)).
 		$$
 \end{observation}
 
 In the following example, we show, that there are semi-simple retractions $r:L\rightarrow L'$,  where $L$ and $L'$ are Lelek fans in $F$ such that $L'\subseteq L$, that are not simple, for which there is a retraction $R:H_L\rightarrow H_{L'}$ such that for each $x\in H_L$, 
 		$$
 		q(R(x)) = r(q(x)).
 		$$
\begin{example}\label{juro}
	Let $L_1$ and $J_1$ be Lelek fans with top $v$ such that $L_1\cap J_1=\{v\}$ as seen in Figure \ref{enkica}. Also, let $L=L_1\cup J_1$ and let $L'=J_1$. Then $L$ and $L'$ are also Lelek fans with top $v$. 
	\begin{figure}[h!]
	\centering
		\includegraphics[width=35em]{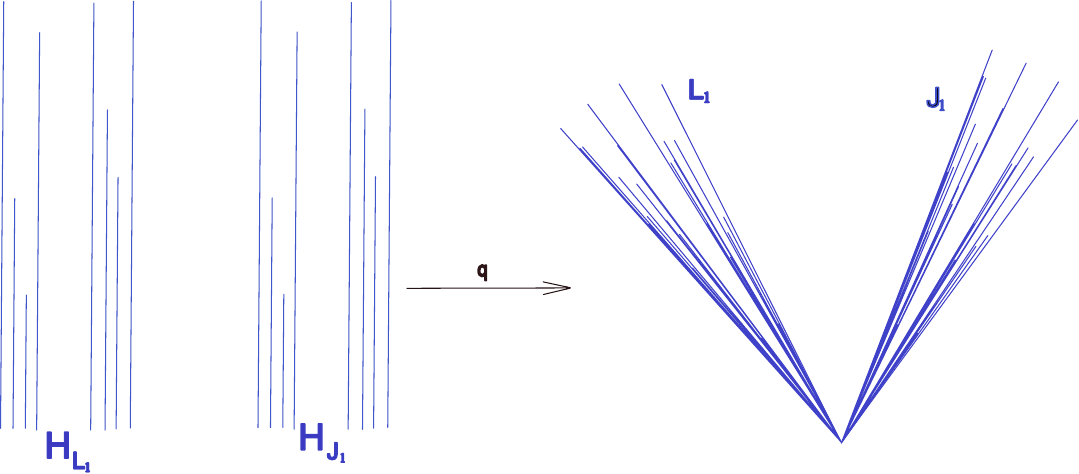}
	\caption{The Lelek fans $L$ and $L'$ from Example \ref{juro}}
	\label{enkica}
\end{figure}
Let $\varphi:L_1\rightarrow J_1$ be a homeomorphism and let $r:L\rightarrow L'$ be defined by 
$$
r(x)=\begin{cases}
				\varphi(x)\text{;} & x\in L_1\\
				x\text{;} & x\in J_1
			\end{cases}
$$
for each $x\in L$. Then $r$ is a semi-simple retraction which is not simple. Observe that $\varphi$ defines a homeomorphism $h:H_{L_1}\rightarrow H_{J_1}$ such that for each $x\in H_{L_1}$, $q(h(x))=\varphi(q(x))$. Note that the function $R:H_L\rightarrow H_{L'}$, defined by 
$$
R(x)=\begin{cases}
				h(x)\text{;} & x\in H_{L_1}\\
				x\text{;} & x\in H_{J_1},
\end{cases}
$$
for each $x\in H_L$, is a retraction such that for each $x\in H_L$, 
 		$$
 		q(R(x)) = r(q(x)).
 		$$
\end{example}
In Example \ref{juro2}, we produce an example of a semi-simple retraction $r:L\rightarrow L'$,  where $L$ and $L'$ are Lelek fans in $F$ such that $L'\subseteq L$, for which there are no retractions $R:H_L\rightarrow H_{L'}$ such that for each $x\in H_L$, 
 		$$
 		q(R(x)) = r(q(x)).
 		$$
\begin{example}\label{juro2}
	For each positive integer $n$, let $L_n$, $J_n$ and $K_n$ be lelek fans, as seen in the Figure \ref{enkica1}, such that
	  \begin{figure}[h!]
	\centering
		\includegraphics[width=35em]{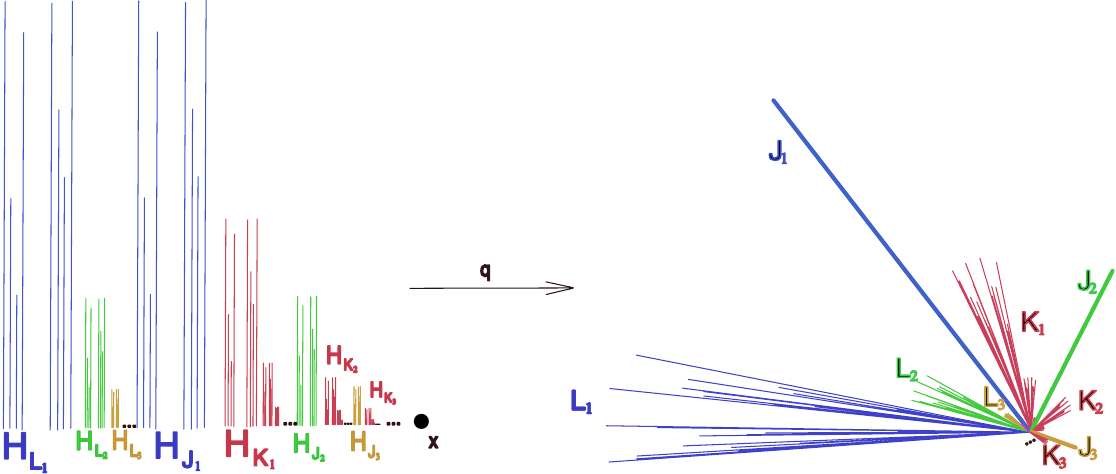}
	\caption{The Lelek fans $L$ and $L'$ from Example \ref{juro2}}
	\label{enkica1}
\end{figure}
 $\displaystyle\lim_{n\to\infty}H_{L_n}$ is a point in $H_{J_1}$, $\displaystyle\lim_{n\to\infty}H_{K_n}=\{x\}$ and $\displaystyle\lim_{n\to\infty}H_{J_n}=\{x\}$. In Figure \ref{enkica1}, the $J_n$'s are pictured as arcs but they are not arcs, they are Lelek fans - imagine that they are subsets of the planes that are perpendicular to the plane, in which $L_1$ is pictured. Also, let 
	$$
	L=\left(\bigcup_{n=1}^{\infty}L_n\right)\cup\left(\bigcup_{n=1}^{\infty}J_n\right)\cup\left(\bigcup_{n=1}^{\infty}K_n\right) ~~~ \textup{ and } ~~~ L'=\left(\bigcup_{n=1}^{\infty}L_n\right)\cup\left(\bigcup_{n=1}^{\infty}K_n\right)
	$$
	 be Lelek fans with top $v$ as seen in Figure \ref{enkica1}. For each positive integer $n$, let $\varphi_n:J_n\rightarrow L_n$ be a homeomorphism and let $r:L\rightarrow L'$ be defined by 
$$
r(x)=\begin{cases}
				\varphi_n(x)\text{;} & \textup{there is a positive integer } n \textup{ such that } x\in J_n\\
				x\text{;} & x\in L'
			\end{cases}
$$
for each $x\in L$. Then $r$ is a semi-simple retraction which is not simple. Note that for each positive integer $n$, $\varphi_n$ defines a homeomorphism $h_n:H_{J_n}\rightarrow H_{L_n}$ such that for each $x\in H_{J_n}$, $q(h_n(x))=\varphi_n(q(x))$. Suppose that there is a retraction $R:H_L\rightarrow H_{L'}$ such that for each $x\in H_L$, 
 		$$
 		q(R(x)) = r(q(x)).
 		$$
 		Then for each positive integer $n$, if $x\in H_{}J_n$, then $R(x)=h_n(x)$. Let $R$ be such a retraction. For each positive integer $n$, let $x_n\in H_{J_n}$ and $y_n\in H_{K_n}$. Then 
 		$$
 		\lim_{n\to\infty}x_n=x ~~~ \textup{ and } ~~~  \lim_{n\to\infty}y_n=x
 		$$ 
 		while 
 		$$
 		\lim_{n\to\infty}R(x_n)\in H_{J_1} ~~~ \textup{ and } ~~~  \lim_{n\to\infty}R(y_n)=\lim_{n\to\infty}y_n=x.
 		$$
 		Since $x\not\in H_{J_1}$, $R$ is not continuous, which is a contradiction. 
\end{example}
The following theorem is our last main result.
 \begin{theorem}\label{sinca}
 	Let $R:H_{J}\rightarrow H_{K}$ be a retraction. Then there is a semi-simple retraction $r:J\to K$ such that for each $x\in H_J$, 
 		$$
 		q(R(x)) = r(q(x)).
 		$$
  \end{theorem}
 \begin{proof}
 Let $\mathcal P=\{\{y\}\times [0,1] \mid y\in Y\}$. Note that for each $P\in \mathcal P$, there is $Q\in \mathcal P$ such that $R(P)\subseteq Q$. For each $[(s,t)]\in J$,  we define
 $$
 r([{(s,t)}])=\begin{cases}
				[R{(s,t)}]\text{;} & [{(s,t)}]\neq v\\
			v\text{;} & [{(s,t)}] = v.
			\end{cases}
 $$
 Then $r$ is a semi-simple retraction from $J$ onto $K$ such that for each ${(s,t)}\in H_J$, 
 		$$
 		q(R{(s,t)}) = r(q{(s,t)}).
 		$$
 \end{proof}
\begin{corollary}
	Let $L$ and $L'$ be Lelek fans in $F$ such that $L'\subseteq L$ and let $R:H_{L}\rightarrow H_{L'}$ be a retraction. Then there is a semi-simple retraction $r:L\to L'$ such that for each $x\in H_L$, 
 		$$
 		q(R(x)) = r(q(x)).
 		$$
 		\end{corollary}
\begin{proof}
	Let $J=L$, $K=L'$ and use Theorem \ref{sinca}.
\end{proof}
 \begin{example}
 	Let $f:[0,1]\rightarrow [0,1]$ be a continuous function such that
 	$$
 	\{(x,y)\in H_L \mid y\leq f(x) \}=H_{L'}
 	$$ 
 	and for each $(x,y)\in H_L$, let
 	$$
 	R(x,y)=\begin{cases}
				(x,y)\text{;} & y\leq f(x)\\
			(x,f(x))\text{;} & y > f(x).
			\end{cases}
 	$$
 	Then $R$ is a retraction from $H_L$ to $H_{L'}$. By Theorem \ref{sinca}, there is a semi-simple retraction from $L$ onto $L'$. In fact, it is a simple retraction from $L$ onto $L'$.
 \end{example}
  \begin{example}
 	Let $W_1=([0,\frac{1}{2}]\times [0,1])\cap H_L$ and let $W_2=([\frac{1}{2},1]\times [0,1])\cap H_L$ and suppose that $W_1\neq \emptyset$, $W_2\neq \emptyset$ and that $\varphi:W_2\to W_1$ is a homeomorphism such that for all $(s,t)\in W_2$, 
 	$$
 	t=0 ~~~ \Longrightarrow ~~~ \varphi(s,t)=(p,0)
 	$$
 	for some $p\in [0,\frac{1}{2}]$. For each $(x,y)\in H_L$, we define
 	$$
 	R(x,y)=\begin{cases}
				(x,y)\text{;} & (x,y)\in W_1\\
			\varphi(x,y)\text{;} & (x,y)\in W_2.
			\end{cases}
 	$$
 	Note that $R$ is a retraction from $H_L$ to $H_{L'}$. By Theorem \ref{sinca}, there is a semi-simple retraction from $L$ onto $L'$. 
 \end{example}
 We conclude the paper with the following open problem.
\begin{problem}
	Find necessary and sufficient conditions for a semi-simple retraction $r:L\rightarrow L'$, for which a retraction $R:H_L\rightarrow H_{L'}$ such that for each $x\in H_L$, 
 		$$
 		q(R(x)) = r(q(x))
 		$$
exists.
\end{problem}

\section{Acknowledgement}
Sina Greenwood is supported by the Marsden Fund Council from Government funding, administered by the Royal Society of New Zealand. This work is also supported in part by the Slovenian Research Agency (research projects J1-4632, BI-HR/23-24-011, and research program P1-0285). 
	

\noindent I. Bani\v c\\
              (1) Faculty of Natural Sciences and Mathematics, University of Maribor, Koro\v{s}ka 160, SI-2000 Maribor,
   Slovenia; \\(2) Institute of Mathematics, Physics and Mechanics, Jadranska 19, SI-1000 Ljubljana, 
   Slovenia; \\(3) Andrej Maru\v si\v c Institute, University of Primorska, Muzejski trg 2, SI-6000 Koper,
   Slovenia\\
             {iztok.banic@um.si}           
     
				\-
				
		\noindent G.  Erceg\\
             Faculty of Science, University of Split, Rudera Bo\v skovi\' ca 33, Split,  Croatia\\
{{goran.erceg@pmfst.hr}       }    

     \-
			
		\noindent S.  Greenwood\\
            Department of Mathematics, University of Auckland, 38 Princes Street, 1010 Auckland, New Zealand\\
{{s.greenwood@auckland.ac.nz}       }

                 	\-
					
  \noindent J.  Kennedy\\
             Department of Mathematics,  Lamar University, 200 Lucas Building, P.O. Box 10047, Beaumont, Texas 77710 USA\\
{{kennedy9905@gmail.com}       }  
\-

	\-
				
		%



\end{document}